\pdfoutput=1
\RequirePackage{ifpdf}
\ifpdf % We are running pdfTeX in pdf mode
\documentclass[pdftex]{sigma}
\else
\documentclass{sigma}
\fi

\usepackage{mathrsfs}
\usepackage{tensor}
\usepackage[vcentermath]{youngtab}

\newcommand{\scaleboxfactor}{0.8}

\newcommand{\alphabox}{\scalebox{\scaleboxfactor}{\ensuremath{\alpha}}}
\newcommand{\betabox}{\scalebox{\scaleboxfactor}{\ensuremath{\beta}}}
\newcommand{\gammabox}{\scalebox{\scaleboxfactor}{\ensuremath{\gamma}}}
\newcommand{\anti}{\young(\alphabox,\betabox,\gammabox)}

\begin{document}

\allowdisplaybreaks

\renewcommand{\thefootnote}{$\star$}

\newcommand{\arXivNumber}{1502.07516}

\renewcommand{\PaperNumber}{024}

\FirstPageHeading

\ShortArticleName{Nijenhuis Integrability for Killing Tensors}

\ArticleName{Nijenhuis Integrability for Killing Tensors\footnote{This paper is a~contribution to the Special Issue
on Analytical Mechanics and Dif\/ferential Geometry in honour of Sergio Benenti.
The full collection is available at \href{http://www.emis.de/journals/SIGMA/Benenti.html}{http://www.emis.de/journals/SIGMA/Benenti.html}}}

\Author{Konrad SCH\"OBEL}
\AuthorNameForHeading{K.~Sch\"obel}

\Address{Mathematisches Institut, Fakult\"at f\"ur Mathematik und Informatik,\\ Friedrich-Schiller-Universit\"at Jena, 07737 Jena, Germany}
\Email{\href{mailto:konrad.schoebel@uni-jena.de}{konrad.schoebel@uni-jena.de}}

\ArticleDates{Received October 30, 2015, in f\/inal form February 26, 2016; Published online March 07, 2016}	

\Abstract{The fundamental tool in the classif\/ication of orthogonal coordinate
	systems in which the Hamilton--Jacobi and other prominent equations can be
	solved by a separation of variables are second order Killing tensors which
	satisfy the Nijenhuis integrability conditions. The latter are a system
	of three non-linear partial dif\/ferential equations. We give a simple and
	completely algebraic proof that for a Killing tensor the third and most
	complicated of these equations is redundant. This considerably simplif\/ies
	the classif\/ication of orthogonal separation coordinates on arbitrary
	(pseudo-)Riemannian manifolds.}


\Classification{70H06;	53A60;	53B20}

\renewcommand{\thefootnote}{\arabic{footnote}}
\setcounter{footnote}{0}

It is a natural problem to classify all coordinate systems in which a given
partial dif\/ferential equation can be solved by a separation of variables~--
the so called \emph{separation coordinates}. For fundamental equations such
as the Hamilton--Jacobi and the Schr\"odinger equation the theory of separation
of variables is built on a characterisation of orthogonal separation
coordinates by second order Killing tensors, i.e., solutions of the Killing
equation
\begin{gather}
	\label{eq:Killing}
	\nabla_\alpha K_{\beta \gamma}+
	\nabla_\beta K_{\gamma\alpha}+
	\nabla_\gamma K_{\alpha\beta }=0,
\end{gather}
which are integrable in the sense that they have simple eigenvalues and the
orthogonal complements of each eigenvector f\/ield form an integrable
distribution. The relation between orthogonal separation coordinates and
integrable Killing tensors was observed in~1891 by Paul St\"ackel in his
Habilitation thesis \cite{Staeckel} and then used by Luther P.~Eisenhart in~1934 to classify orthogonal se\-paration coordinates on 3-dimensional Euclidean
space and on the 3-dimensional sphere~\cite{Eisenhart}. His results were
generalised to arbitrary spaces of constant curvature of any dimension by
Kalnins and Miller in~1986~\cite{Kalnins, Kalnins&Miller}.

The property of an endomorphism f\/ield to be integrable in the above sense has
been cast into the form of a system of three non-linear partial dif\/ferential
equations by Nijenhuis in 1950~\cite{Nijenhuis}. Explicitly, an endomorphism
$K$ with simple eigenvalues is integrable if and only if it satisf\/ies the
\emph{Nijenhuis integrability conditions}
\begin{subequations}
	\label{eq:Nijenhuis}
	\begin{gather}
		0 =N\indices{^\delta_{[\beta\gamma}}g_{\alpha]\delta},\\
		0 =N\indices{^\delta_{[\beta\gamma}}K_{\alpha]\delta},\\
		0 =N\indices{^\delta_{[\beta\gamma}}K_{\alpha]\varepsilon}K\indices{^\varepsilon_\delta},
	\end{gather}
\end{subequations}
where the square brackets stand for antisymmetrisation and $N$ denotes the
\emph{Nijenhuis torsion} of~$K$, def\/ined by
\begin{gather*}
	N(X,Y) := K^2[X,Y]-K[KX,Y]-K[X,KY]+[KX,KY]
\end{gather*}
and given in local coordinates by{\samepage
\begin{gather}
 \label{eq:torsion}
	N\indices{^\alpha_{\beta\gamma}}
	=K\indices{^\alpha_\delta}K\indices{^\delta_{[\beta;\gamma]}}
	+K\indices{^\delta_{[\beta}}K\indices{^\alpha_{\gamma];\delta}},
\end{gather}
where a semicolon denotes a covariant derivative.}

Of course, these equations have not been known to St\"ackel or Eisenhart.
Neither did they play any role in the complete classif\/ication of Kalnins and
Miller. However, they reveal that the classif\/ication of orthogonal separation
coordinates is actually an \emph{algebraic geometric problem}. Indeed, the
Nijenhuis integrability conditions \eqref{eq:Nijenhuis} are algebraic in~$K$
and its derivatives and hence constitute a set of homogeneous algebraic
equations on the space of second order Killing tensors. Note that, as the
space of solutions to the overdetermined linear equation~\eqref{eq:Killing},
this space is a f\/inite-dimensional vector space. Consequently, the set of
Killing tensors satisfying the Nijenhuis integrability conditions is a
projective variety. Moreover, this variety comes along with a natural group
action, namely the isometry group of the manifold.

To be more precise, orthogonal separation coordinates are in one-to-one
correspondence with so-called \emph{St\"ackel systems}, i.e., $n$-dimensional
spaces of integrable Killing tensors which mutually commute in the algebraic
sense. This leads to the following remarkable observation~\cite{Schoebel16}:
\begin{quote}
	For any (pseudo-)Riemannian manifold the set of orthogonal separation
	coordinates carries a very rich structure: It is a projective variety
	isomorphic to a subvariety in the Grassmannian of $n$-dimensional
	subspaces in the space of Killing tensors, equipped with a natural action
	of the isometry group.
\end{quote}
To the best of our knowledge, this point has never been made in the literature
and the structure of these varieties had never been elucidated. The reason is
probably that a general solution of the equations~\eqref{eq:Nijenhuis} was
deemed intractable~\cite{Horwood&McLenaghan&Smirnov}.

\looseness=-1
Recently it has been possible to rewrite the Nijenhuis integrability
conditions explicitly as algebraic equations for constant curvature manifolds~\cite{Schoebel12} and to solve them in the least non-trivial case, namely for
the sphere of dimension three~\cite{Schoebel14}. The outcome, a detailed
algebraic geometric description of the variety of integrable Killing tensors
as well as the variety of St\"ackel systems, has lead to a surprising connection
between separation coordinates on spheres on one hand and algebraic curves on
the other. More precisely, the set of orthogonal separation coordinates
modulo isometries on the $n$-dimensional sphere is naturally parametrised by
the real Deligne--Mumford moduli space $\bar{\mathscr M}_{0,n+2}(\mathbb R)$ of
stable algebraic curves of genus $0$ with $n+2$ marked points~\cite{Schoebel&Veselov}. As a consequence, the well known classif\/ication of
Kalnins and Miller can be reinterpreted in terms of the geometry and
combinatorics of Stashef\/f polytopes. This also revealed a~hitherto unknown
operad structure on equivalence classes of separation coordinates on spheres.
To date, comparable results are unknown for manifolds other than spheres.

An important step in the explicit solution of the Nijenhuis equations has been
the proof that the third of the Nijenhuis conditions is redundant when applied
to a Killing tensor. According to a footnote in
\cite{Horwood&McLenaghan&Smirnov}, this had previously been proven by Steve
Czapor for Euclidean space in dimension three using Gr\"obner bases. The author
extended this to arbitrary constant curvature manifolds~\cite{Schoebel12}.
The purpose of the present paper is to give a simple proof that this result
holds in full generality. This will considerably simplify the classif\/ication
of orthogonal separation coordinates on arbitrary manifolds.

\begin{theorem}
	For a Killing tensor on an arbitrary Riemannian manifold the third of the
	Nijenhuis equations \eqref{eq:Nijenhuis} is redundant. The same holds
	true on a pseudo-Riemannian manifold.
\end{theorem}

\begin{remark}
	For a Killing tensor the f\/irst two Nijenhuis equations are in general
	independent~\cite{Schoebel14}.
\end{remark}

\begin{remark}
	Note that a St\"ackel system contains Killing tensors whose eigenvalues are
	not simple (the metric for instance). This is why we did not impose
	simple eigenvalues in the above theorem.
\end{remark}

\begin{remark}
	Instead of the Nijenhuis conditions \eqref{eq:Nijenhuis}, the vanishing of
	the Haantjes tensor
	\begin{gather*}
		H(X,Y):=K^2N(X,Y)-KN(KX,Y)-KN(X,KY)+N(KX,KY)
	\end{gather*}
	is often used as a criterion for integrability. Being of order four in
	$K$, this condition is of the same complexity as the third of the
	Nijenhuis equations, while the f\/irst two are only of order two and three.
	Our result therefore shows that for Killing tensors the Nijenhuis
	conditions are better suited than the Haantjes tensor.
\end{remark}

\section*{The proof}

\looseness=-1
We will prove the statement pointwise. That is, we f\/ix an arbitrary point~$p$
in the manifold and consider the restrictions $K_{\alpha\beta}(p)$ and
$(\nabla_\alpha K_{\beta\gamma})(p)$ of the Killing tensor f\/ield and its
covariant derivative to the tangent space at~$p$. The statement then becomes
a~purely algebraic statement on these two tensors. For ease of notation we
omit to indicate the dependence on the chosen point~$p$.

A Killing tensor is symmetric by def\/inition. Hence at the f\/ixed point we can
choose an orthonormal basis of the tangent space in which the Killing tensor
is diagonal, i.e., ${K^\alpha}_\beta=\lambda_\alpha{\delta^\alpha}_\beta$ (no
sum). In this basis, the Nijenhuis torsion \eqref{eq:torsion} reads
\begin{gather*}
	N_{\alpha\beta\gamma}
	=\young(\betabox,\gammabox)(\lambda_\alpha-\lambda_\gamma)K_{\alpha\beta;\gamma},
\end{gather*}
where the Young projector $\young(\betabox,\gammabox)$ stands for
antisymmetrisation in $\beta$ and $\gamma$. Substituted into the Nijenhuis
integrability conditions \eqref{eq:Nijenhuis} we get
\begin{gather*}
	0
	=N\indices{^\delta_{[\beta\gamma}}g_{\alpha]\delta}
	=\anti N_{\alpha\beta\gamma}
	=\anti(\lambda_\alpha-\lambda_\gamma)K_{\alpha\beta;\gamma},\\
	0
	=N\indices{^\delta_{[\beta\gamma}}K_{\alpha]\delta}
	=\anti\lambda_\alpha N_{\alpha\beta\gamma}
	=\anti\lambda_\alpha(\lambda_\alpha-\lambda_\gamma)K_{\alpha\beta;\gamma},\\
	0
	=N\indices{^\delta_{[\beta\gamma}}K_{\alpha]\varepsilon}K\indices{^\varepsilon_\delta}
	=\anti\lambda_\alpha^2N_{\alpha\beta\gamma}
	=\anti\lambda_\alpha^2(\lambda_\alpha-\lambda_\gamma)K_{\alpha\beta;\gamma},
\end{gather*}
where the Young projector
\begin{gather*}
	\anti
\end{gather*}
denotes a complete antisymmetrisation in the indices $\alpha$, $\beta$ and~$\gamma$. Using that $K_{\alpha\beta;\gamma}$ is symmetric in $\alpha$, $\beta$
and that a complete antisymmetrisation over $\alpha$, $\beta$, $\gamma$ can be split
into an antisymmetrisation in $\alpha$, $\beta$ and a subsequent sum over the
cyclic permutations of~$\alpha$, $\beta$, $\gamma$, we can rewrite the preceding
equations as
\begin{subequations}
	\label{eq:Ki}
	\begin{gather}
		\label{eq:K1}0=(\lambda_\alpha-\lambda_\beta)K_{\alpha\beta;\gamma}+\text{cyclic},\\
		\label{eq:K2}0=(\lambda_\alpha-\lambda_\beta)(\lambda_\alpha+\lambda_\beta-\lambda_\gamma)K_{\alpha\beta;\gamma}+\text{cyclic},\\
		\label{eq:K3}0=(\lambda_\alpha-\lambda_\beta)\bigl((\lambda_\alpha+\lambda_\beta)^2-\lambda_\alpha\lambda_\beta
-\lambda_\beta\lambda_\gamma-\lambda_\gamma\lambda_\alpha\bigr)K_{\alpha\beta;\gamma}+\text{cyclic},
	\end{gather}
\end{subequations}
where ``$+$~cyclic'' stands for a summation over the cyclic permutations
of~$\alpha$, $\beta$, $\gamma$. These equations are one by one equivalent to the
Nijenhuis integrability conditions~\eqref{eq:Nijenhuis}. In the same way we
can write the Killing equation as
\begin{gather}
	\label{eq:K0}
	0=K_{\alpha\beta;\gamma}+\text{cyclic}.
\end{gather}
Adding appropriate multiples of \eqref{eq:K1} to \eqref{eq:K2} and
\eqref{eq:K3}, we can simplify \eqref{eq:Ki} to
\begin{subequations}
	\begin{gather}
		\label{eq:K1'}0 =(\lambda_\alpha-\lambda_\beta)K_{\alpha\beta;\gamma}+\text{cyclic},\\
		\label{eq:K2'}0 =(\lambda_\alpha-\lambda_\beta)(\lambda_\alpha+\lambda_\beta)K_{\alpha\beta;\gamma}+\text{cyclic},\\
		\label{eq:K3'}0 =(\lambda_\alpha-\lambda_\beta)(\lambda_\alpha+\lambda_\beta)^2K_{\alpha\beta;\gamma}+\text{cyclic}.
	\end{gather}
\end{subequations}
We want to prove that \eqref{eq:K0} together with \eqref{eq:K1} and
\eqref{eq:K2} imply \eqref{eq:K3}, which is equivalent to prove that
\eqref{eq:K0} together with \eqref{eq:K1'} and \eqref{eq:K2'} imply
\eqref{eq:K3'}. To this end we write the f\/irst three equations in matrix form
as
\begin{gather*}
	\begin{bmatrix}
		1&1&1\\
		\lambda_\alpha -\lambda_\beta &\lambda_\beta -\lambda_\gamma &\lambda_\gamma- \lambda_\alpha \\
		\lambda_\alpha^2-\lambda_\beta^2&\lambda_\beta^2-\lambda_\gamma^2&\lambda_\gamma^2-\lambda_\alpha^2
	\end{bmatrix}
	\begin{bmatrix}
		K_{\alpha\beta ;\gamma}\\
		K_{\beta \gamma;\alpha}\\
		K_{\gamma\alpha;\beta }
	\end{bmatrix}
	=0.
\end{gather*}
The determinant of the coef\/f\/icient matrix is an antisymmetric cubic polynomial
in $\lambda_\alpha$, $\lambda_\beta$, $\lambda_\gamma$ and hence a multiple of the
Vandermode determinant. If the eigenvalues
$\lambda_\alpha$, $\lambda_\beta$, $\lambda_\gamma$ are pairwise dif\/ferent, this
implies that
$K_{\alpha\beta;\gamma}=K_{\beta\gamma;\alpha}=K_{\gamma\alpha;\beta}=0$. If
exactly two of the eigenvalues are equal, say
$\lambda_\alpha\not=\lambda_\beta=\lambda_\gamma$, then we have
$K_{\alpha\beta;\gamma}=-\frac12K_{\beta\gamma;\alpha}=K_{\gamma\alpha;\beta}$.
For three equal eigenvalues the only restriction on $K_{\alpha\beta;\gamma}$
is the Killing equation~\eqref{eq:K0}. In all three cases we see that the
equation~\eqref{eq:K3'} is also satisf\/ied.

For a pseudo-Riemannian manifold the statement follows from the above and the
fact that in the space of symmetric tensors the (complex) diagonalisable ones
are dense. Indeed, the above reasoning for
$C_{\alpha\beta}=K_{\alpha\beta}(p)$ remains true even if~$C_{\alpha\beta}$ is
(complex) diagonalisable, but not necessarily the restriction of a Killing
tensor to the tangent space at~$p$. By continuity, it is therefore also true
if $C_{\alpha\beta}$ is not diagonalisable and therefore also for the
restriction of an arbitrary Killing tensor.

\subsection*{Acknowledgements}
The author would like to acknowledge the anonymous referees for their
contribution to improve the paper.

\pdfbookmark[1]{References}{ref}
\LastPageEnding


\begin{thebibliography}{99}
\footnotesize\itemsep=0pt

\bibitem{Eisenhart}
Eisenhart L.P., Separable systems of {S}t\"ackel, \href{http://dx.doi.org/10.2307/1968433}{\textit{Ann. of Math.}}
 \textbf{35} (1934), 284--305.

\bibitem{Horwood&McLenaghan&Smirnov}
Horwood J.T., McLenaghan R.G., Smirnov R.G., Invariant classif\/ication of
 orthogonally separable Hamiltonian systems in Euclidean space, \href{http://dx.doi.org/10.1007/s00220-005-1331-8}{\textit{Comm.
 Math. Phys.}} \textbf{259} (2005), 670--709, \href{http://arxiv.org/abs/math-ph/0605023}{math-ph/0605023}.

\bibitem{Kalnins}
Kalnins E.G., Separation of variables for {R}iemannian spaces of constant
 curvature, \textit{Pitman Monographs and Surveys in Pure and Applied
 Mathematics}, Vol.~28, Longman Scientif\/ic \& Technical, Harlow, 1986.

\bibitem{Kalnins&Miller}
Kalnins E.G., Miller Jr. W., Separation of variables on {$n$}-dimensional
 {R}iemannian manifolds. {I}.~{T}he {$n$}-sphe\-re~{$S_n$} and {E}uclidean
 {$n$}-space~{${\bf R}^n$}, \href{http://dx.doi.org/10.1063/1.527088}{\textit{J.~Math. Phys.}} \textbf{27} (1986),
 1721--1736.

\bibitem{Nijenhuis}
Nijenhuis A., {$X_{n-1}$}-forming sets of eigenvectors, \href{http://dx.doi.org/10.1016/S1385-7258(51)50028-8}{\textit{Indag. Math.}}
 \textbf{54} (1951), 200--212.

\bibitem{Schoebel12}
Sch{\"o}bel K., Algebraic integrability conditions for {K}illing tensors on
 constant sectional curvature manifolds, \href{http://dx.doi.org/10.1016/j.geomphys.2012.01.006}{\textit{J.~Geom. Phys.}} \textbf{62}
 (2012), 1013--1037, \href{http://arxiv.org/abs/1004.2872}{arXiv:1004.2872}.

\bibitem{Schoebel14}
Sch{\"o}bel K., The variety of integrable {K}illing tensors on the 3-sphere,
 \href{http://dx.doi.org/10.3842/SIGMA.2014.080}{\textit{SIGMA}} \textbf{10} (2014), 080, 48~pages, \href{http://arxiv.org/abs/1205.6227}{arXiv:1205.6227}.

\bibitem{Schoebel16}
Sch{\"o}bel K., Are orthogonal separable coordinates really classif\/ied?,
\textit{SIGMA} \textbf{12} (2016), to appear, \href{http://arxiv.org/abs/1510.09028}{arXiv:1510.09028}.

\bibitem{Schoebel&Veselov}
Sch{\"o}bel K., Veselov A.P., Separation coordinates, moduli spaces and
 {S}tashef\/f polytopes, \href{http://dx.doi.org/10.1007/s00220-015-2332-x}{\textit{Comm. Math. Phys.}} \textbf{337} (2015),
 1255--1274, \href{http://arxiv.org/abs/1307.6132}{arXiv:1307.6132}.

\bibitem{Staeckel}
St{\"a}ckel P.,
Die {I}ntegration der {H}amilton--{J}acobischen {D}ifferentialgleichung mittelst {S}eparation der {V}a\-riab\-len,
Habilitationsschrift, Universit\"at Halle, 1891.

\end{thebibliography}
\end{document}